\documentclass[11pt,reqno]{amsart}

\usepackage{hyperref}
\usepackage{graphicx}
\usepackage{color}
\usepackage{amsfonts,amssymb}
\usepackage{bbm} 

\usepackage{mathrsfs}

\usepackage{a4wide}

\newtheorem{neu}{}[section]
\newtheorem{Cor}[neu]{Corollary}
\newtheorem*{Cor*}{Corollary}
\newtheorem{Thm}[neu]{Theorem}
\newtheorem*{Thm*}{Theorem}
\newtheorem{Prop}[neu]{Proposition}
\newtheorem*{Prop*}{Proposition}
\theoremstyle{definition}
\newtheorem{Lemma}[neu]{Lemma}
\newtheorem*{Rmk*}{Remark}
\newtheorem{Rmk}[neu]{Remark}

\newtheorem*{Ex*}{Example}

\newtheorem*{Qu*}{Question}

\newtheorem{Def}[neu]{Definition}

\newcommand{\Z}{\mathbb{Z}}
\newcommand{\R}{\mathbb{R}}

\newcommand{\pf}{\longrightarrow}

\newcommand{\id}{\mathrm{id}}

\newcommand{\om}{\omega}

\newcommand{\A}{\mathcal{A}}
\newcommand{\At}{\widetilde{\mathcal{A}}}

\newcommand{\F}{\mathcal{F}}
\newcommand{\Ft}{\widetilde{\mathcal{F}}}

\newcommand{\M}{\mathcal{M}}

\newcommand{\E}{\mathcal{E}}
\renewcommand{\L}{\mathscr{L}}

\renewcommand{\H}{\mathrm{H}}

\newcommand{\Ham}{\mathrm{Ham}}

\newcommand{\CF}{\mathrm{CF}}

\newcommand{\HF}{\mathrm{HF}}

\newcommand{\HFl}{\mathrm{HF}_\mathrm{loc}}
\newcommand{\RFH}{\mathrm{RFH}}
\newcommand{\CFl}{\CF_{\mathrm{loc}}}

\newcommand{\Crit}{\mathrm{Crit}}

\newcommand{\beq}{\begin{equation}}
\newcommand{\beqn}{\begin{equation}\nonumber}
\newcommand{\eeq}{\end{equation}}

\newcommand{\bea}{\begin{equation}\begin{aligned}}
\newcommand{\bean}{\begin{equation}\begin{aligned}\nonumber}
\newcommand{\eea}{\end{aligned}\end{equation}}

\numberwithin{equation}{section}

\definecolor{Urs}{rgb}{0,.7,0}
\definecolor{Peter}{rgb}{0,0,1}
\definecolor{red}{rgb}{1,0,0}

\newcommand{\He}{\mathscr{H}}

\begin{document}
\title{Leaf-wise intersections and Rabinowitz Floer homology}
\author{Peter Albers}
\author{Urs Frauenfelder}
\address{
    Peter Albers\\
    Department of Mathematics\\
    Purdue University}
\email{palbers@math.purdue.edu}
\address{
    Urs Frauenfelder\\
    Department of Mathematics and Research Institute of Mathematics\\
    Seoul National University}
\email{frauenf@snu.ac.kr}
\begin{abstract}
In this article we explain how critical points of a particular perturbation of the Rabinowitz action functional give rise to leaf-wise intersection points in hypersurfaces of restricted contact type. This is used to derive existence and multiplicity results for leaf-wise intersection points in hypersurfaces of restricted contact type in general exact symplectic manifolds. The notion of leaf-wise intersection points was introduced by Moser \cite{Moser_A_fixed_point_theorem_in_symplectic_geometry}.
\end{abstract}
\maketitle

\section{Introduction}

We consider a closed hypersurface $\Sigma\subset(M,\om=d\lambda)$ in an exact symplectic manifold $(M,\om)$ such that $(\Sigma,\alpha:=\lambda_{|\Sigma})$ is a contact manifold. Moreover, we assume that $\Sigma$ bounds a compact region in $M$ and that $M$ is convex at infinity, that is, $M$ is isomorphic to the symplectization of a compact contact manifold at infinity. $\Sigma$ is foliated by the leaves of the characteristic line bundle which is spanned by the Reeb vector field $R$ of $\alpha$. For $x\in\Sigma$ we denote by $L_x$ the leaf through $x$. Furthermore, we denote by $\Ham_c(M,\om)$ the group of compactly supported Hamiltonian diffeomorphism. The following question was addressed by Moser \cite{Moser_A_fixed_point_theorem_in_symplectic_geometry}.
\begin{Qu*}\label{question:existence of leaf-wise intersection}
Given  $\phi\in\Ham_c(M,\om)$, does there exist a leaf-wise intersection point, that is, $x\in\Sigma$ with $\phi(x)\in L_x$?
\end{Qu*}

\begin{Def}
We denote by $\wp(\Sigma,\alpha)>0$ the minimal period of a Reeb orbit of $(\Sigma,\alpha)$ which is contractible in $M$. If there exists no contractible Reeb orbit we set $\wp(\Sigma,\alpha)=\infty$. 
\end{Def}

Our first main result is the following.\\[-1.5ex]

\noindent\textbf{Theorem A.}  \textit{If $\phi\in\Ham_c(M,\om)$ has Hofer norm $||\phi||<\wp(\Sigma,\alpha)$, then there exists a leaf-wise intersection point for $\phi$.}

\begin{Rmk}
The case $(M,\Sigma)=(\R^{2n},S^{2n-1})$ shows that Theorems A is sharp since $\wp(S^{2n-1})$ equals the displacement energy of the sphere $S^{2n-1}$. In particular, the smallness assumption in Theorem A is necessary.
\end{Rmk}

The proof of Theorem A uses a stretching of the neck argument for gradient flow lines of a perturbed, time dependent Rabinowitz action functional. More sophisticatedly, using local Rabinowitz Floer homology around the action value 0 we obtain the following multiplicity result.\\[-1.5ex]

\noindent\textbf{Theorem B.}  \textit{For a generic Hamiltonian diffeomorphism $\phi\in\Ham_c(M,\om)$ with $||\phi||<\wp(\Sigma,\alpha)$
\beq
\#\{\text{leaf-wise intersection points}\}\geq\sum_i b_i(\Sigma,\Z/2)\;.
\eeq
}\\[-1ex]

If the full Rabinowitz Floer homology is non-zero we obtain much stronger results.  For the construction of Rabinowitz Floer homology we refer the reader to \cite{Cieliebak_Frauenfelder_Restrictions_to_displaceable_exact_contact_embeddings}, see also Section \ref{sec:RabinowitzFH}.\\[-1ex]

\noindent\textbf{Theorem C.} \textit{If the Rabinowitz Floer homology of $(M,\Sigma)$ does not vanish, $\RFH(M,\Sigma)\neq0$, then there always exists a leaf-wise intersection point for $\psi\in\Ham_c(M,\om)$.\\[.5ex]
}\\[-1ex]
We point out that we make no assumption on the Hofer norm of $\psi$. Moreover, as mentioned above, without the assumption $\RFH(M,\Sigma)\neq0$ Theorem C does not hold in general.

\begin{Rmk}
In the article \cite{Cieliebak_Frauenfelder_Restrictions_to_displaceable_exact_contact_embeddings} examples with non-vanishing Rabinowitz Floer homology are provided. See \cite{Cieliebak_Frauenfelder_Oancea_Rabinowitz_Floer_homology_and_symplectic_homology} for further examples.
\end{Rmk}

\begin{Rmk}
The leaf-wise intersection points found in Theorems A and B are always contractible in the following sense. For any Hamiltonian function $H:S^1\times M\pf\R$ such that $\phi=\phi_H$ the leaf-wise intersection point can be completed to a loop $\gamma$ by first following the flow of $\phi_H^t$ and then the Reeb flow in such a way that $\gamma$ is contractible in $M$, see Lemma \ref{lemma:leafwise intersection is contractible} below.
\end{Rmk}

\begin{Rmk}
As in Theorem B local Rabinowitz Floer homology around the action value of a non-contractible Reeb orbit can be considered. Similar techniques then lead to existence results for non-contractible leaf-wise intersections points. In fact, generically each Reeb orbit gives rise to two different leaf-wise intersection points since the local homology is isomorphic to the homology of a circle.
\end{Rmk}

\subsection{History of the problem and related results}

The problem addressed above is a special case of the leaf-wise coisotropic intersection problem. For that let $N\subset(M,\om)$ be a coisotropic submanifold. Then $N$ is foliated by isotropic leafs. The problem asks for a leaf $L$ such that $\phi(L)\cap L\neq\emptyset$ for $\phi\in\Ham_c(M,\om)$. 

The first result was obtained by Moser in \cite{Moser_A_fixed_point_theorem_in_symplectic_geometry} for simply connected $M$ and $C^1$-small $\phi$. This was later generalized by Banyaga \cite{Banyaga_On_fixed_points_of_symplectic_maps} to non-simply connected $M$.

The $C^1$-smallness assumption was replaced by Hofer, Ekeland-Hofer in \cite{Hofer_On_the_topological_properties_of_symplectic_maps},\cite{Ekeland_Hofer_Two_symplectic_fixed_point_theorems_with_applications_to_Hamiltonian_dynamics} for hypersurfaces of restricted contact type in $\R^{2n}$ by a much weaker smallness assumption, namely that the Hofer norm of $\phi$ is smaller than a certain symplectic capacity. Only recently, the result by Ekeland-Hofer was generalized in two different directions. It was extended by Dragnev \cite{Dragnev_Symplectic_rigidity_symplectic_fixed_points_and_global_perturbations_of_Hamiltonian_systems} to so-called ``coisotropic submanifolds of contact type in $\R^{2n}$''.  Among other results Ginzburg \cite{Ginzburg_Coisotropic_Intersections} generalized from restricted contact type in $\R^{2n}$ to restricted contact type in subcritical Stein manifolds. Moreover, examples by Ginzburg \cite{Ginzburg_Coisotropic_Intersections} show that the Ekeland-Hofer result is a symplectic rigidity result, namely it becomes wrong for arbitrary hypersurfaces. Recently Ziltener \cite{Ziltener_coisotropic} and Gurel \cite{Gurel_leafwise_coisotropic_intersection} obtained results on leaf-wise intersection points using entirely different methods from this article.

Theorem A gives a complete answer to the existence problem of leaf-wise intersection points within the class of restricted contact type hypersurfaces. After this article was published on the arXiv an independent proof of Theorem A was given by Gurel \cite{Gurel_leafwise_coisotropic_intersection}.

Theorem B is (to the authors' knowledge) the first time that a general multiplicity result for leaf-wise intersection points in the Hofer-small case is proved. In the $C^1$-small case multiplicity results were obtained by Moser and Banyaga. In the special case of fibrations Ziltener \cite{Ziltener_coisotropic} proves similar multiplicity results to Theorem B.

Theorem C is (again to the authors' knowledge) the first global (i.e.~valid for all Hamiltonian diffeomorphisms) existence result for leaf-wise intersection points. 

\subsubsection*{Acknowledgments}
We thank Felix Schlenk for helpful comments.
This article was written during visits of the first author at the Ludwig-Maximilians-Universit\"at M\"unchen and visits of the second author at the ETH Z\"urich. Both authors thank these institutions for their stimulating working atmospheres. The authors are partially supported by the German Research Foundation (DFG) through Priority Program 1154 "Global Differential Geometry", grant FR 2637/1-1, and NSF grant DMS-0805085.

\section{A perturbation of the Rabinowitz action functional}\label{sec:RFH}

We recall that $\Sigma\subset(M,\om=d\lambda)$ is a closed hypersurface in an exact symplectic manifold such that $(\Sigma,\alpha=\lambda_{|\Sigma})$ is a contact manifold. Moreover, $\Sigma$ is assumed to bound a compact region in $M$. We denote by $R$ the Reeb vector field of $\alpha$. Moreover, we define the vector field $Y$ by $d\lambda(Y,\cdot)=\lambda(\cdot)$. 

\begin{Lemma}\label{lemma:def liouville vector field}
The vector field $Y$ is a Liouville vector field for $(\Sigma,\alpha)$, that is, $\mathcal{L}_Y\om=\om$ and $Y\pitchfork\Sigma$. In particular, $(\Sigma,\alpha)$ is of restricted contact type.
\end{Lemma}

\begin{proof}
Since $\lambda(Y)=d\lambda(Y,Y)=0$ we compute $\mathcal{L}_Y\lambda=d(\iota_Y\lambda)+\iota_Yd\lambda=\lambda$. Since the Lie derivative commutes with the exterior differential we conclude $\mathcal{L}_Y\om=\om$. If we assume that $Y\in T_x\Sigma$ then $0=d\alpha(Y,R)=d\lambda(Y,R)=\lambda(R)=\alpha(R)=1$. This contradiction shows $Y\pitchfork\Sigma$.
\end{proof}

The flow $\phi_Y^t$ of the Liouville vector field is defined near $\Sigma$. We fix $\delta_0>0$ so that ${\phi_Y^t}|_{\Sigma}$ is defined for all $|t|\leq\delta_0$ and define a function $\widehat{G}$ by $\widehat{G}(\phi_Y^t(x))=t$ for all $x\in\Sigma$. For $0<\delta\leq\delta_0$ we set 
\beq
U_\delta:=\{x\in M\mid |\widehat{G}(x)|<\delta\}.
\eeq
Since $\Sigma$ bounds we can choose a $G:M\pf\R$ which is locally constant outside $U_{\delta_0}$, coincides with $\widehat{G}$ on $U_{\frac{\delta_0}{2}}$, and such that $G^{-1}(0)=\Sigma$. Thus, the Hamiltonian vector field $X_G$ of $G$ satisfies ${X_G}|_{\Sigma}=R$. Finally, we fix a smooth function $\rho:S^1\pf\R$ with $\int_0^1\rho(t)dt=1$ and $\mathrm{supp}(\rho)\subset(0,\tfrac12)$ and set
\beq\label{eqn:def_of_F}
F(t,x):=\rho(t)G(x)\;.
\eeq
Therefore, the Hamiltonian vector fields satisfy
\beq
X_F(t,x)=\rho(t)X_G(x)\;.
\eeq 
We recall the definition of the positive and negative part of the Hofer norm.
\begin{Def}
Let $H:S^1\times M\pf\R$ a compactly supported Hamiltonian function. We set
\beq
||H||_+:=\int_0^1\max_{x\in M} H(t,x) dt\qquad||H||_-:=-\int_0^1\min_{x\in M} H(t,x) dt=||-H||_+
\eeq
and 
\beq
||H||=||H||_++||H||_-\;.
\eeq
For $\phi\in\Ham_c(M,\om)$ the Hofer norm is 
\beq
||\phi||=\inf\{||H||\mid \phi=\phi_H\}\;.
\eeq
\end{Def}

\begin{Lemma}\label{lemma:norms equivalent}
For all $\phi\in\Ham_c(M,\om)$ 
\beq
||\phi||=|||\phi|||:=\inf\{||H||\mid \phi=\phi_H,\;H(t,\cdot)=0\;\;\forall t\in[0,\tfrac12]\}\;.
\eeq
\end{Lemma}

\begin{proof}
$||\phi||\leq|||\phi|||$ is obvious. To prove the reverse inequality pick a smooth monotone map $r:[0,1]\to[0,1]$ with $r(\tfrac12)=0$ and $r(1)=1$. For $H$ with $\phi_H=\phi$ we set $H^r(t,x):=r'(t)H(r(t),x)$. Then a direct computation shows $\phi_{H^r}=\phi_H$, $||H^r||=||H||$, and $H^r(t,x)=0$ for all $t\in[0,\tfrac12]$. This finishes the proof.
\end{proof}

From now on we assume that $H(t,\cdot)=0$ for all $t\in[0,\tfrac12]$. Then for $F$ as in equation \eqref{eqn:def_of_F} the perturbed Rabinowitz action functional is defined as follows 
\beq\label{eqn:Rabinowitz action functional}
\A_H^F(u,\eta):=-\int_0^1u^*\lambda-\int_0^1H(t,u(t))dt-\eta\int_0^1F(t,u(t))dt
\eeq
where $u\in C^{\infty}(S^1,M)$ and $\eta\in\R$. Critical points $(u,\eta)\in\Crit\A_H^F$ satisfy 
\beq\label{eqn:critical point equation}\left.
\begin{aligned}
&\partial_tu=X_H(t,u)+\eta X_F(t,u)\\[1ex]
&\int_0^1F(t,u)dt=0
\end{aligned}
\;\;\right\}
\eeq
In the following proposition we observe that existence of a critical point of $\A_H^F$ gives rise to a positive answer of the leaf-wise intersection problem mentioned in the introduction.

\begin{Prop}\label{prop:critical point answers question}
Let $(u,\eta)\in\Crit\A_H^F$. Then $x=u(\tfrac12)$ satisfies $\phi_H(x)\in L_x$. Thus, $x$ is a leaf-wise intersection point.
\end{Prop}

\begin{proof}
For $t\in[0,\tfrac12]$ we compute, using $H(t,\cdot)=0$ for all $t\leq\tfrac12$,
\bea
\frac{d}{dt}G(u(t))&=dG(u(t))\cdot\partial_tu\\
&=dG(u(t))\cdot[\underbrace{X_H(t,u)}_{=0}+\eta \underbrace{X_F(t,u)}_{=\rho(t)X_G(u)}]\\
&=0\;,
\eea
since $dG(X_G)=0$. Hence $G(u(t))=c=$const for $t\leq\tfrac12$. Thus,
\bea
0=\int_0^1F(t,u)dt=\int_0^1\rho(t)G(u(t))dt=c
\eea
Therefore, $G(u(t))=c=0$, and since $G^{-1}(0)=\Sigma$ we have  $u(t)\in\Sigma$ for $t\in[0,\tfrac12]$. In particular, $u(\tfrac12),u(0)=u(1)\in\Sigma$.

For $t\in[\tfrac12,1]$ we have $F(t,\cdot)=0$. Thus, the loop $u$ solves the equation $\partial_tu=X_H(t,u)$ on $[\tfrac12,1]$, and therefore, $u(1)=\phi_H(u(\frac12))$. We conclude that $\phi_H(u(\frac12))\in\Sigma$.
Using again that for $t\in[0,\tfrac12]$, $\partial_tu=X_H(t,u)+\eta X_F(t,u)=\eta X_F(t,u)=\eta \rho(t) X_G(u)$ and $u(t)\in\Sigma$ we see that $u(1)=u(0)\in L_{u(\frac12)}$ since $X_G|_\Sigma=R$. 

With the definition $x:=u(\frac12)$ we then have $\phi_H(x)=u(1)\in L_x$. This concludes the proof.
\end{proof}

In the following we establish necessary analytical properties of the perturbed Rabinowitz action functional. For later purposes we allow from now on the function $H$ to be $s$-dependent as follows: $H_s(t,x)=H_-(t,x)$ for $s\leq-1$ and $H_s(t,x)=H_+(t,x)$ for $s\geq1$. Moreover, $H_s(t,\cdot)=0$ for $t\in[0,\tfrac12]$, and $H_s$ has compact support uniformly in $s$. We choose a family $J(s,t)$ of compatible almost complex structures on $M$ such that $J(s,t)$ is independent of $s$ for $s\leq-1$ and $s\geq1$. The norm of the gradient of $\A_H^F$ equals
\beq\label{eqn:norm_of_gradient}
||\nabla\A_H^F(u,\eta)||^2=||\partial_tu-X_{H_s}(t,u)-\eta X_F(t,u)||_{L^2}^2+\Big|\int_0^1F(t,u(t))dt\Big|^2
\eeq
where the $L^2$ norm is taken with respect to the metric $g_{(s,t)}(\cdot,\cdot):=\om(\cdot,J(s,t)\cdot)$. We denote by $\L$ the component of the contractible loops in $M$.

\begin{Def}\label{def:gradient flow line}
A gradient flow line of $\A_H^F$ is (formally) a map $w=(u,\eta)\in C^{\infty}(\R,\L\times\R)$
solving the ODE 
\beq\label{eqn:gradient flow line}
\partial_s w(s)+ \nabla_s\A_H^F(w(s))=0\;,
\eeq
where the gradient is taken with respect to metric $\mathfrak{m}_s$ defined as follows. Let $(\hat{u}_1,\hat{\eta}_1)$ and $(\hat{u}_2,\hat{\eta}_2)$ be two tangent vectors in $T_{(u,\eta)}(\L\times\R)$. We set
\beq
\mathfrak{m}_s\big((\hat{u}_1,\hat{\eta}_1),\,(\hat{u}_2,\hat{\eta}_2)\big):=\int_0^1g_{(s,t)}\big(\hat{u}_1,\hat{u}_2\big)dt+\hat{\eta}_1\hat{\eta}_2\;.
\eeq
According to Floer's interpretation, \cite{Floer_The_unregularized_gradient_flow_of_the_symplectic_action}, this means that $u$ and $\eta$ are smooth maps $u:\R\times S^1\pf M$ and $\eta:\R\pf\R$ solving
\beq\label{eqn:gradient flow equation}\left.
\begin{aligned}
&\partial_su+J(s,t,u)\big(\partial_tu-X_{H_s}(t,u)-\eta X_F(t,u)\big)=0\\[1ex]
&\partial_s\eta-\int_0^1F(t,u)dt=0.
\end{aligned}
\;\;\right\}
\eeq
\end{Def}

\begin{Def}
The energy of a map $w\in C^{\infty}(\R,\L\times\R)$ is defined as
\beq
E(w):=\int_{-\infty}^\infty||\partial_s w||^2ds\;.
\eeq
\end{Def}

\begin{Lemma}\label{lemma:energy estimate for gradient lines}
Let $w$ be a gradient flow line of $\nabla_s\A_H^F$. Then 
\beq\label{eqn:energy estimate for gradient lines}
E(w)\leq\A_{H_-}^F(w_-)-\A_{H_+}^F(w_+)+\int_{-\infty}^\infty||\partial_sH_s||_-ds\;.
\eeq 
Moreover, equality holds if $\partial_sH_s=0$. 
\end{Lemma}

\begin{proof}
It follows from the gradient flow equation \eqref{eqn:gradient flow line}
\bea
E(w)&=-\int_{-\infty}^\infty d\A_{H_s}^F(w(s))[\partial_sw]ds\\[1ex]
&=-\int_{-\infty}^\infty \frac{d}{ds}\Big(\A_{H_s}^F(w(s))\Big)ds+\int_{-\infty}^\infty \big(\partial_s\A_{H_s}^F\big)(w)ds\\[1ex]
&=\A_{H_-}^F(w_-)-\A_{H+}^F(w_+)-\int_{-\infty}^\infty\int_0^1\partial_sH_s(t,u)dtds\\[1ex]
&\leq\A_{H_-}^F(w_-)-\A_{H_+}^F(w_+)+\int_{-\infty}^\infty||\partial_sH_s||_-ds\;.
\eea
\end{proof}

\begin{Lemma}\label{lemma:energy estimate for gradient lines II}
Let $w$ be a gradient flow line of $\nabla_s\A_H^F$. Then 
\beq
|\A_H^F(w(s_0))|\leq\max\{\A_{H_-}^F(w_-),-\A_{H_+}^F(w_+)\}+\int_{-\infty}^\infty||\partial_sH_s||_-ds
\eeq 
for all $s_0\in\R$.
\end{Lemma}

\begin{proof}
The proof follows from the proof of Lemma \ref{lemma:energy estimate for gradient lines} by replacing in the first line $\int_{-\infty}^\infty$ by $\int_{-\infty}^{s_0}$ resp.~$\int_{s_0}^\infty$, and $E(w)\geq0$.
\end{proof}

\begin{Thm}\label{thm:action bounds imply compactness}
Let $w_n=(u_n,\eta_n)$ be a sequence of gradient flow lines for which there exists $a<b$ such that
\beq
a\leq\A_H^F\big(w_n(s)\big)\leq b\qquad\forall s\in\R\;.
\eeq
Then for every reparametrisation sequence $\sigma_n\in\R$ the sequence $w_n(\cdot+\sigma_n)$ has a subsequence which converges in $C^\infty_\mathrm{loc}(\R,\L\times\R)$.
\end{Thm}

\begin{proof}
The proof follows from standard arguments in Floer theory as soon as we establish
\begin{enumerate}
\item a uniform $L^\infty$ bound on $u_n$,
\item a uniform $L^\infty$ bound on $\eta_n$,
\item a uniform $L^\infty$ bound on the derivatives of $u_n$.
\end{enumerate} 
Indeed, assuming (1)--(3) bootstrapping the gradient-flow equation will establish $C^\infty_\mathrm{loc}$-convergence, see \cite[Appendix B.4]{McDuff_Salamon_J_holomorphic_curves_and_symplectic_topology}. The $L^\infty$ bound on $u_n$ follows from the convexity at infinity of $(M,\om)$. Once the $L^\infty$ bound on $\eta_n$ has been established, the $L^\infty$ bound on the derivatives of $u_n$ follows in the following way. If the derivatives would explode we would obtain non-constant holomorphic spheres as limits, see \cite[Chapter 4.2]{McDuff_Salamon_J_holomorphic_curves_and_symplectic_topology}. But in an exact symplectic manifold non-constant holomorphic spheres don't exist. The $L^\infty$ bound on $\eta_n$ is the content of the following proposition.
\end{proof}

\begin{Prop}\label{prop:L infinity bounds for eta}
Given critical points $w_-,w_+\in\Crit\A_H^F$ there exists a constant $\kappa=\kappa(w_-,w_+)$ such that every gradient flow line $w=(u,\eta)$ of $\A_H^F$ with $\lim_{s\to\pm\infty}=w_\pm$ satisfies
\beq
||\eta||_{L^\infty(\R)}\leq \kappa\;.
\eeq
\end{Prop}

The proof of the proposition goes along the same lines as in \cite[Corollary 3.3]{Cieliebak_Frauenfelder_Restrictions_to_displaceable_exact_contact_embeddings} and relies on the following lemma. The proof of the proposition is given after the proof of the lemma.

\begin{Lemma}\label{lemma:crucial lemma for eta bound}
There exists $\epsilon>0$ and $C>0$ such that for all $(u,\eta)\in C^{\infty}(S^1,M)\times\R$ we have
\beq
||\nabla_s\A_H^F(u,\eta)||<\epsilon\quad\Longrightarrow\quad|\eta|\leq C\Big(|\A_H^F(u,\eta)|+1\Big)
\eeq
where the norm of the gradient is given in equation \eqref{eqn:norm_of_gradient}.
\end{Lemma}

\begin{proof}
We will use again the notation introduced below the proof of Lemma \ref{lemma:def liouville vector field}. We fix $0<2\delta<\min\{1,\delta_0\}$, in particular, we have $G(x)=\widehat{G}(x)$ for $x\in U_\delta$.

\underline{Claim 1}: Assume that $u(t)\in U_\delta$ for all $t\in(\tfrac12,1)$. Then there exists a constant $C_1>0$
\beq
|\eta|\leq C_1\Big(|\A_H^F(u,\eta)|+||\nabla_s\A_H^F(u,\eta)||+1\Big)\;.
\eeq

\begin{proof}[Proof of Claim 1]
We compute
\bean
|\A_H^F(u,\eta)|&=\left|-\int_0^1u^*\lambda-\int_0^1H(s,t,u(t))dt-\eta\int_0^1F(t,u(t))dt\right|\\[1ex]
&\geq\left|-\int_0^1u^*\lambda\right|-||H||_{L^\infty}-\delta|\eta|\\[1ex]
&\geq-||H||_{L^\infty}-\delta|\eta|+\left|\int_0^1\lambda(u(t))\big[\partial_tu-X_{H_s}(t,u)-\eta X_F(t,u)\big]dt\right.\\
&\qquad+\int_{\tfrac12}^1\lambda(u(t))\big[X_{H_s}(t,u)\big]dt+\underbrace{\int_0^{\tfrac12}\lambda(u(t))\big[\eta \underbrace{X_F(t,u)}_{=\rho(t)R(u(t))}\big]dt}_{=\eta}\bigg|\\
&\geq|\eta|-\delta|\eta|-C_{\lambda,\delta}\;||\partial_tu-X_{H_s}(t,u)-\eta X_F(t,u)||_{L^1}-C_{\lambda, H}\\[1ex]
&\geq\tfrac12|\eta|-C_{\lambda,\delta}\;||\partial_tu-X_{H_s}(t,u)-\eta X_F(t,u)||_{L^2}-C_{\lambda, H}\\[1ex]
&\geq\tfrac12|\eta|-C_{\lambda,\delta}\;||\nabla_s\A_H^F(u,\eta)||_{L^2}-C_{\lambda, H}
\eea
where $C_{\lambda,\delta}:=||\lambda_{|U_\delta}||_{L^\infty}$ and $C_{\lambda, H}:=||H||_{L^\infty}+C_{\lambda,\delta}||X_{H}||_{L^\infty}$. This inequality immediately implies Claim 1.
\end{proof}

\underline{Claim 2}: There exists $\epsilon=\epsilon(\delta)$ with the following property. If for $(u,\eta)$ there exists $t\in[0,\tfrac12]$ with $G(u(t))\geq\delta$ then $||\nabla_s\A_H^F(u,\eta)||\geq\epsilon$

\begin{proof}[Proof of Claim 2]
If in addition $G(u(t))\geq\frac{\delta}{2}$ holds for all $t\in[0,\tfrac12]$ then using \eqref{eqn:norm_of_gradient}
\beq
||\nabla_s\A_H^F(u,\eta)||\geq\Big|\int_0^1F(t,u(t))dt\Big|\geq\frac{\delta}{2}\int_0^1\rho(t)dt=\frac{\delta}{2}
\eeq
since $F(t,x)=\rho(t)G(x)$. Otherwise there exists $t'\in[0,\tfrac12]$ with $G(u(t'))<\frac{\delta}{2}$. Thus, we can find $0\leq a<b\leq\tfrac12 $ such that for all $t\in[a,b]$
\beq
\frac{\delta}{2}=G(u(a))\leq G(u(t))\leq G(u(b))=\delta
\eeq
or 
\beq
\delta=G(u(a))\geq G(u(t))\geq G(u(b))=\frac{\delta}{2}\;.
\eeq
We only treat the first case. The second is completely analogous.
\bea
||\nabla_s\A_H^F(u,\eta)||&\geq||\partial_tu-X_{H_s}(t,u)-\eta X_F(t,u)||_{L^2}\\[1ex]
&\geq\left(\int_a^b||\partial_tu-\underbrace{X_{H_s}(t,u)}_{=0}-\eta X_F(t,u)||^2dt\right)^{\tfrac12}\\[1ex]
&\geq\left(\int_a^b\frac{1}{||\nabla G||^2}\big|g_t(\partial_tu,\nabla G)-\eta \underbrace{g_t(X_F(t,u),\nabla G)}_{=0}\big|^2dt\right)^{\tfrac12}\\[1ex] 
&\geq\frac{1}{||\nabla G||_{L^\infty}}\left(\int_a^b\left|\frac{d}{dt}G(u(t))\right|^2dt\right)^{\tfrac12}\\[1ex]
&\geq\frac{1}{||\nabla G||_{L^\infty}}\int_a^b\left|\frac{d}{dt}G(u(t))\right|dt\\[1ex]
&\geq\frac{1}{||\nabla G||_{L^\infty}}\int_a^b\frac{d}{dt}G(u(t))dt\\[1ex]
&=\frac{\delta}{2||\nabla G||_{L^\infty}}
\eea
where we used $g_t(X_F,\nabla G)=dG(X_F)=\rho(t)dG(X_G)=\rho(t)\om(X_G,X_G)=0$. Since $||\nabla G||_{L^\infty}\geq1$ we set $\epsilon(\delta):=\frac\delta2$. This proves Claim 2.
\end{proof}
Setting $3\delta:=\min\{1,\delta_0\}$, $\epsilon:=\epsilon(\delta)$ according to Claim 2, and $C:=C_1+\epsilon$, $C_1$ as in Claim 1, the lemma follows.
\end{proof}

\begin{proof}{(of Proposition \ref{prop:L infinity bounds for eta})}\quad
From Lemma \ref{lemma:energy estimate for gradient lines} it follows that
\beq
E(w)\leq\A_{H_-}^F(w_-)-\A_{H_+}^F(w_+)+C_H
\eeq
where $C_H:=\int||\partial_sH_s||_-ds$. We fix $\epsilon$ and $C$ as in Lemma \ref{lemma:crucial lemma for eta bound}. For $\sigma\in\R$ we define
\beq
\tau(\sigma):=\inf\{\tau\geq0\mid||\nabla_s\A_H^F(w(\sigma+\tau))||\leq\epsilon\}
\eeq
and compute
\beq
E(w)\geq\int_{\sigma}^{\sigma+\tau(\sigma)}\underbrace{||\nabla_s\A_H^F(w(s))||^2}_{\geq\epsilon^2}ds\geq\tau(\sigma)\epsilon^2
\eeq
From the second equation in \eqref{eqn:gradient flow equation} it follows directly that
\beq
||\partial_s\eta||_{L^\infty}\leq||F||_{L^\infty}\;.
\eeq
The norm $||F||_{L^\infty}$ is finite since $dF=0$ outside a compact set. Finally according to Lemma \ref{lemma:energy estimate for gradient lines II}
\beq
|\A_H^F(w(s))|\leq\max\big\{\A_H^F(w_-),-\A_H^F(w_+)\big\}+C_H=:\Xi
\eeq
The last three inequalities together with Lemma \ref{lemma:crucial lemma for eta bound} imply
\bea
|\eta(\sigma)|&\leq|\eta(\sigma+\tau(\sigma))|+\int_{\sigma}^{\sigma+\tau(\sigma)}|\partial_s\eta|ds\\[1ex]
&\leq C\Big(\left|\A_H^F\big(w(\sigma+\tau(\sigma))\big)\right|+1\Big)+||F||_{L^\infty}\tau(\sigma)\\[1ex]
&\leq C(\Xi+1)+||F||_{L^\infty}\frac{E}{\epsilon^2}\;.
\eea
This proves the proposition.
\end{proof}

\subsection{Proof of Theorem A}

Recall that $\A:\E\pf\R$ is a functional and $C\subset\Crit{\A}$ then $C$ is called a Morse-Bott component if the following two conditions hold.

\begin{itemize}
 \item $C$ is a submanifold of $\E$ 
 \item For all $c\in C$ we have $T_cC=\ker\He_\A(c)$ where $\He_\A$ is the Hessian of $\A$.
\end{itemize}

\begin{Lemma}\label{lemma:Morse_Bott}
The subset $\Sigma\subset\Crit{\A^F}$ is a Morse-Bott component. 
\end{Lemma}

\begin{proof}
Let $c=(p,0)$ with $p\in\Sigma$. An element $(\hat{v},\hat{\eta})\in C^{\infty}(S^1,T_pM)\times\R$ is in the kernel of the Hessian $\He_\A^F(c)$ if and only if it solves the following equations
\beq\left.\begin{aligned}
\partial_t\hat{v}=\hat{\eta}\rho(t)X_G(p)\\
\int_0^1\rho(t)dG(p)\hat{v}dt=0
\end{aligned}\right\}
\eeq 
Integrating the first equation we obtain
\beq
\hat{v}(1)=\hat{v}(0)+\hat{\eta}X_G(p)\;.
\eeq
Using that $\hat{v}$ is a loop and $X_G(p)\neq0$ we conclude $\hat{\eta}=0$ and thus $\hat{v}=\hat{v}_0\in T_pM$. The second equation implies
\beq
dG(p)\hat{v}_0=0
\eeq
and therefore $\hat{v}_0\in T_p\Sigma=\ker dG(p)$.
\end{proof}

\begin{proof}
We choose $H:S^1\times M\pf\R$ such that $\phi_H=\phi$, $||H||<\wp(\Sigma,\alpha)$, and such that $H(t,x)=0$ for all $t\in[0,\tfrac12]$, see Lemma \ref{lemma:norms equivalent}. For $r\geq0$ we choose a smooth family of functions $\beta_r\in C^\infty(\R,[0,1])$ satisfying
\begin{enumerate}
\item for $r\geq1$: $\beta_r'(s)\cdot s\leq0$ for all $s\in\R$, $\beta_r(s)=1$ for $|s|\leq r-1$, and $\beta_r(s)=0$ for $|s|\geq r$, 
\item for $r\leq1$: $\beta_r(s)\leq r$ for all $s\in\R$ and $\mathrm{supp}\beta_r\subset[-1,1]$,
\item $\lim_{r\to\infty}\beta_r(s\mp r)=:\beta_\infty^\pm(s)$ exists, where the limit is taken with respect to the $C^\infty$ topology.
\end{enumerate} 
We set 
\beq\label{eqn:homotopy_from_proof_of_ThmA}
K_r(s,t,x):=\beta_r(s)H(t,x)\;.
\eeq
We fix a point $p\in\Sigma$ and consider the moduli space 
\beq
\M:=\left\{(r,w)\in[0,\infty)\times C^{\infty}(\R,\L\times\R)\bigg|\; \begin{aligned}&w\text{ solves \eqref{eqn:gradient flow equation} for $K_r$}\\&\lim_{s\to-\infty}w(s)=(p,0),\;\lim_{s\to\infty}w(s)\in\Sigma \end{aligned}\right\}\;.
\eeq


\underline{Claim}: If there exists no leaf-wise intersection point, then $\M$ is compact. Moreover, its boundary consists of the point $(0,p,0)$ only.\\

Assuming the Claim we prove the assertion of Theorem A. The moduli space $\M$ is the zero-set of a Fredholm section of a Banach-space bundle over a Banach manifold. Its index equals 1. Moreover, since by Lemma \ref{lemma:Morse_Bott} $\A^F=\A^{K_0}$ is Morse-Bott along $\Sigma$ the Fredholm section is regular at the boundary point $(0,p,0)$. It is well-known that a Fredholm section can be perturbed to a transverse Fredholm section given that its zero-set is compact. Since the Fredholm section is already transverse at the boundary point it suffices to perturb away from the boundary. Thus, assuming the claim we obtain from $\M$ a smooth compact manifold with boundary being the point $(0,p,0)$. Such a manifold does not exists. Thus, to finish the proof it remains to show the claim. 

According to Lemma \ref{lemma:energy estimate for gradient lines} we have for $(r,w)\in\M$ with $\lim_{s\to\infty}w=(p',0)$
\bea
E(w)&\leq\A^F_0(p,0)-\A^F_0(p',0)+\int_{-\infty}^\infty||\partial_sK_r||_-ds\\[1ex]
&=\int_{-\infty}^\infty||\beta_r'(s)H||_-ds\\[1ex]
&=\int_{-\infty}^0\beta_r'(s)||H||_-ds-\int_{0}^\infty\beta_r'(s)||H||_+ds\\[1ex]
&=\beta_r(0)\big(||H||_-+||H||_+\big)\\[1ex]
&\leq||H||\;.
\eea
Moreover, if $r=0$ in the above then $E(w)=\A^F_0(p,0)-\A^F_0(p',0)=0$ thus $w$ is constant and therefore $w(s)=(p,0)$ for all $s\in\R$. In particular, for $r=0$ the only solution in $\M$ is $w(s)=(p,0)$. Finally, since $\A^F_0(p,0)=\A^F_0(p',0)=0$ the above computation implies 
\beqn
-||H||\leq\A^{K_r}(w(s))\leq||H|| \quad\forall s\in\R\;.
\eeq
Since we have uniform action bounds we know by Theorem \ref{thm:action bounds imply compactness} that the sequence $w_n$ converges (after choosing a subsequence) to some solution $v$ of the gradient flow equation. In case that  $v(+\infty)\not\in\Sigma$  (the other case is analogous) we choose an open set $U\subset \L \times \R$ containing only the constant critical points. Let $s_n$ be the first time so that $w_n(\sigma_n)\not \in U$, i.e.~the first exit time. This is well-defined for large enough $n$ since $v(+\infty)\not\in\Sigma$. Now consider the reparametrised sequence $u_n := \sigma_n * w_n := w_n(\cdot + \sigma_n)$. By $C^\infty_{\mathrm{loc}}$ compactness the sequence $u_n$ converges to $u$ (after choice of a subsequence). Then $u$ is a non-constant gradient flowline since $u(0) \in \partial U$ and $u(-\infty)=(p,0)$ using again that $\Sigma$ is Morse-Bott, see Lemma \ref{lemma:Morse_Bott}. Thus, one of the following has to exist
\begin{enumerate}
\item a non-constant gradient flow line $v$ of $\A^F_0$ with one asymptotic end being $(p,0)$,
\item a gradient flow line $v$ of $\A_{\beta_{\pm}^{\infty}H}^F$, where $\beta_{\pm}^{\infty}$ is as above.
\end{enumerate}
Moreover, $E(v)\leq\limsup E(w_n)\leq||H||<\wp(\Sigma,\alpha)$. In the first case the sequence $r_n$ converges, whereas in the second case $r_n\to\infty$.
If there exists no leaf-wise intersection points, then the second case cannot occur since otherwise one asymptotic end of $v$ is a critical point of $\A_H^F$ which gives a leaf-wise intersection point according to Proposition \ref{prop:critical point answers question}.

In the first case not both asymptotic ends of $v$ can be of the form $(q,0)\in\Crit\A_0^F$  for some $q\in\Sigma$ since otherwise $E(v)=0$ according to Lemma \ref{lemma:energy estimate for gradient lines}. Hence the gradient flow line would be constant. Therefore, one asymptotic end of $v$ is of the form $(\gamma,\eta)$ where $\gamma$ is a Reeb orbit (contractible in $M$) of period $\eta\neq0$. Thus, $|\eta|=|\A^F_0(\gamma,\eta)|=E(v)<\wp(\Sigma,\alpha)$, which contradicts the definition of $\wp(\Sigma,\alpha)$. This finishes the proof.
\end{proof}

\begin{Lemma}\label{lemma:leafwise intersection is contractible}
For any Hamiltonian function $H:S^1\times M\pf\R$ such that $\phi=\phi_H$ the leaf-wise intersection point found in Theorem A can be completed to a loop $\gamma$ by first following the flow of $\phi_H^t$ and then the Reeb flow in such a way that $\gamma$ is contractible in $M$.
\end{Lemma}

\begin{proof}
From the previous proof it follows immediately that for the chosen Hamiltonian function $H$ the leaf-wise intersection point $x=u(0)$, where $(u,\eta)\in\Crit\A_H^F$, can be completed to a contractible loop $\gamma$. We observe that $\gamma(t)=\phi_{H+\eta F}^t(x)$ is a 1-periodic orbit of the Hamiltonian function $H+\eta F$. If $K$ is another Hamiltonian function with $\phi=\phi_K$, then $H+\eta F$ and $K+\eta F$ have the same time-1-maps $\phi_{H+\eta F}=\phi_{K+\eta F}$. Moreover, it follows from the existence of at least one contractible periodic orbit for the time-1-map of the flow $\phi^{-t}_{H+\eta F}\circ\phi^t_{K+\eta F}$ that the loop $\phi_{H+\eta F}^t(x)$ is contractible if and only if $\phi_{K+\eta F}^t(x)$ is contractible. The existence of a contractible periodic orbit for $\phi^{-t}_{H+\eta F}\circ\phi^t_{K+\eta F}$ follows from Floer's proof of the Arnold conjecture \cite{Floer_Morse_theory_for_Lagrangian_intersections}.

\end{proof}

\subsection{The perturbed Rabinowitz action functional is generically Morse}

We set 
\beq\label{eqn:def_of_fancy_H}
\mathcal{H}:=\{H\in C^\infty(S^1\times M)\mid H(t,\cdot)=0\;\forall t\in[0,\tfrac12]\}\;.
\eeq
The specific time support of functions $H\in\mathcal{H}$ is crucial in the proof of Proposition \ref{prop:critical point answers question}. Fortunately, the set $\mathcal{H}$ still generates $\Ham(M,\om)$, see Lemma \ref{lemma:norms equivalent}. Moreover, it is large enough so that the perturbed Rabinowitz action functional is generically Morse.
\begin{Thm}\label{thm:Morse}
For a generic $H\in\mathcal{H}$ the perturbed Rabinowitz action functional $\A_H^F$ is Morse.
\end{Thm}

\begin{proof}
The proof is postponed to the appendix \ref{appendix}.
\end{proof}

\begin{Rmk}
It is straight forward to prove that $\A_H^F$ is Morse if one does not insist that $H$ has time support in $[\tfrac12,1]$, see \cite{Cieliebak_Frauenfelder_Paternain_Symplectic_Topology}. The proof of the genericity of the Morse property follows a standard scheme once it is shown that a certain linear operator is surjective. This operator is composed out of two summands. One is the Hessian of $\A_H^F$ and the other comes from the variation in $H$. Without restrictions on the time support surjectivity follows essentially directly from examining the second summand. In the situation of this paper this fails and we crucially use the Hessian.
\end{Rmk}

\subsection{Rabinowitz Floer homology}\label{sec:RabinowitzFH}

The definition of Floer homology $\HF(\A_H^F)$ proceeds as usual. We choose an $s$-independent $H:S^1\times M\pf\R$. In addition, we require that $H(t,\cdot)=0$ for all $t\in[0,\tfrac12]$. Moreover, we assume that the perturbed Rabinowitz action functional $\A_H^F$ is Morse.  Then we define the $\Z/2$ vector space
\beq
\CF(\A_H^F):=\Big\{\xi=\!\!\!\!\sum_{c\in\Crit\A_H^F}\!\!\xi_cc\mid\xi_c\in\Z/2,\;\#\{c\mid\xi_c\neq0, \A_H^F(c)\geq\kappa\}<\infty\quad\forall\kappa\in\R\Big\}
\eeq
and the moduli space
\beq
\M(c_-,c_+):=\Big\{w\mid w\text{ solves the gradient flow equation \eqref{eqn:gradient flow equation}},\;\lim_{s\to\pm\infty}=c_\pm\Big\}/\R\;.
\eeq
Transversality for the moduli spaces $\M(c_-,c_+)$ can be achieved by abstract perturbation theory. For example, one can use the theory of polyfolds by Hofer-Wysocki-Zehnder. In fact, since there are no automorphism groups acting, the framework of $M$-polyfolds \cite{HWZ_A_General_Fredholm_Theory_I,HWZ_A_General_Fredholm_Theory_II} is sufficient to achieve transversality. Indeed, the space of broken trajectories is an $M$-polyfold and the gradient flow equation is a Fredholm section into an $M$-polybundle. The moduli space $\M(c_-,c_+)$ is the zero set of this Fredholm section. Using the abstract perturbation theory developed in \cite{HWZ_A_General_Fredholm_Theory_II} then achieves transversality.

It also is conceivable that $\M(c_-,c_+)$ is smooth for a generic choice of an $S^1$-family $J(t,\cdot)$ of compatible almost complex structures similarly as in the work of \cite{Floer_Hofer_Salamon_Transversality_in_elliptic_Morse_theory_for_the_symplectic_action}.

By abuse of notation the smooth manifold obtained by perturbing the gradient flow equation will again be denoted by $\M(c_-,c_+)$. We set $n(c_-,c_+)$ to be the $\Z/2$-number of elements in the zero-dimensional components of $\M(c_-,c_+)$. Then the linear map defined on generators by
\bea
\partial:\CF(\A_H^F)&\pf\CF(\A_H^F)\\
c&\mapsto \sum_{d}n(c,d)d
\eea
satisfies $\partial^2=0$. We set $\HF(\A_H^F):=\H(\CF(\A_H^F),\partial)$. 

\begin{Thm}\label{thm:perturbed RFH equals RFH}
If $H$ is such that $\A_H^F$ is a Morse function then
\beq
\HF(\A_H^F)\cong\RFH(M,\Sigma)\;.
\eeq
\end{Thm}

\begin{proof}
We choose an $s$-dependent homotopy from $H$ to $0$. Because of Theorem \ref{thm:action bounds imply compactness} the usual continuation homomorphisms are well-defined and isomorphisms. We conclude
\beq
\HF(\A_H^F)\cong\HF(\A_0^F)=\RFH(M,\Sigma)\;.\qedhere
\eeq
\end{proof}

\subsection{Proof of Theorem C}

Theorem C from the introduction is the following corollary of Theorem \ref{thm:perturbed RFH equals RFH}.
\begin{Cor}\label{cor:theorem a'}
If $\RFH(M,\Sigma)\neq0$, then there exists a leaf-wise intersection point for any $\phi\in\Ham_c(M,\om)$.
\end{Cor}

\begin{proof}
We assume by contradiction that there exists no leaf-wise intersection points. In particular, by Proposition \ref{prop:critical point answers question} $\Crit\A_H^F=\emptyset$ and thus $\A_H^F$ is Morse with $0=\HF(\A_H^F)\cong\RFH(M,\Sigma)$. This contradiction finishes the proof.
\end{proof}

\subsection{Local Rabinowitz Floer homology}

In the following we assume that $\phi_H\in\Ham_c(M,\om)$ is chosen so that $\A_H^F$ is Morse. For $||H||<\wp(\Sigma,\alpha)$ we define 
\beqn
\Crit_\mathrm{loc}(\A_H^F):=\Big\{(u,\eta)\in\Crit\A_H^F\mid u \text{ is contractible in }M,\; -||H||_+\leq\A_H^F(u,\eta)\leq||H||_-\Big\}\;.
\eeq
We note that the set $\Crit_\mathrm{loc}(\A_H^F)$ is finite. This follows from the Arzela-Ascoli theorem since the Lagrange multiplier $\eta$ is uniformly bounded according to Lemma \ref{lemma:crucial lemma for eta bound}.  We define the finite dimensional $\Z/2$ vector space
\beq
\CFl(\A_H^F):=\Crit_\mathrm{loc}(\A_H^F)\otimes\Z/2\;.
\eeq
$(\CFl(\A_H^F),\partial_\mathrm{loc})$ is a differential complex since the action along gradient flow lines is decreasing. Define local Rabinowitz Floer homology by $\HFl(\A_H^F):=\H(\CFl(\A_H^F),\partial_\mathrm{loc})$. 
 
%



\begin{Def}
We abbreviate the number of leaf-wise intersection points of $\phi_H\in\Ham_c(M,\om)$ by $\nu_{\mathrm{leaf}}(\phi_H)$.
\end{Def}

\begin{Lemma}\label{lemma:crucial}
If $\A_H^F$ is Morse and $||H||<\wp(\Sigma,\alpha)$ the inequalities
\beq
\nu_{\mathrm{leaf}}(\phi_H)\geq\dim\CFl(\A_H^F)\geq\dim\HFl(\A_H^F)\;.
\eeq
hold.
\end{Lemma}

\begin{proof}
The second inequality is obvious. To prove the first inequality we point out that two critical points $(u,\eta)\neq(u',\eta')\in\Crit(\A_H^F)$ can give rise to the same leaf-wise intersection point only if the underlying leaf of the Reeb flow is closed. Indeed, if $u(\frac12)=u'(\frac12)$ then according to Proposition \ref{prop:critical point answers question} we know that $u(1)=\phi_H(u(\frac12))=\phi_H(u'(\frac12))=u'(1)$. Trivially, $u(1)=u(0)=u'(0)=u'(1)$.  This is illustrated in figure \ref{fig:non-unique}. The map $u$ would be following the solid line, whereas the map $u'$ would follow the solid arc and the dotted part of the Reeb orbit. 

\begin{figure}[htb]
\input{leafwise1.pstex_t}\caption{}\label{fig:non-unique}
\end{figure}

We denote by $(u|_{[0,\frac12]})^-$ the path $u|_{[0,\frac12]}$ traversed in the opposite direction. Then, the map
\beq
\gamma:=u'|_{[0,\frac12]}\#(u|_{[0,\frac12]})^-
\eeq
is a closed loop in $\Sigma$ which (up to reparametrization) is a non-trivial Reeb orbit. The loop $\gamma$ is contractible in $M$ since $\gamma$ is homotopic to the loop $u'\#u^-$ which is the concatenation of two contractible loops and thus contractible. Next, we compute using $u|_{[\frac12,1]}=u'|_{[\frac12,1]}$
\bea
\left|\A_H^F(u,\eta)-\A_H^F(u',\eta')\right|&=\left|-\int_0^\frac12 u^*\alpha-\int_\frac12^1 u^* \lambda-\eta\int_0^\frac12\rho(t)\underbrace{G(u(t))}_{=0}dt-\int_\frac12^1H(t,u(t))dt\right.\\
	&\left.\;-\left(\int_0^\frac12 (u')^*\alpha-\int_\frac12^1 (u')^*\lambda-\eta'\int_0^\frac12\rho(t)\underbrace{G(u'(t))}_{=0}dt-\int_\frac12^1H(t,u'(t))dt\right)\right|\\
	&=\left|\int_\frac12^1 (u')^*\alpha-\int_\frac12^1 u^*\alpha\right|\\
	&=\left|\int_{S^1}\gamma^*\alpha\right|\geq\wp(\Sigma,\alpha)>||H||
\eea
If we assume that $(u,\eta)\neq(u',\eta')\in\Crit_\mathrm{loc}(\A_H^F)$ give rise to the same leaf-wise intersection, then by definition of $\Crit_\mathrm{loc}(\A_H^F)$  we have 
\beq
||H||=||H||_++||H||_-\geq\left|\A_H^F(u,\eta)-\A_H^F(u',\eta')\right|>||H||\;.
\eeq
This contradiction finishes the proof.
\end{proof}

\subsection{Proof of Theorem B}

Theorem B from the introduction follows from Theorem \ref{thm:Morse}, Lemma \ref{lemma:crucial}, and the following proposition.

\begin{Prop}
If $\phi_H\in\Ham_c(M,\om)$ satisfies  $||\phi_H||<\wp(\Sigma,\alpha)$ and if $\A_H^F$ is Morse, then there exists an injective homomorphism 
\beq
\theta:\H(\Sigma;\Z/2)\pf\HFl(\A_H^F)\;.
\eeq
\end{Prop}

\begin{proof}
We first observe that
\beq
\H(\Sigma;\Z/2)\cong\HFl(\A_0^F)\;. 
\eeq
Indeed, this follows from the fact that locally around the action value $0$ the Rabinowitz action functional $\A_0^F$ is Morse-Bott with critical manifold $\Sigma$, on which the action functional vanishes. Since the functional only has one critical value the complex of cascades, see \cite{Frauenfelder_Arnold_Givental_Conjecture}, computing the Morse-Bott homology $\HFl(\A_0^F)$ equals the Morse complex of the critical manifold $\Sigma$. Let  $\theta$ be the continuation homomorphism $\theta:\HFl(\A_0^F)\pf\HFl(\A_H^F)$ in local Floer homology. With formula \eqref{eqn:energy estimate for gradient lines} and $||H||_++||H||_-=||H||<\wp(\Sigma,\alpha)$ one checks via an energy-action estimate that $\theta$ is well-defined when using the homotopy $\beta^+_\infty(s)H$ from the proof of Theorem A, see equation \eqref{eqn:homotopy_from_proof_of_ThmA}.  The same energy-action estimate shows that the reverse continuation homomorphism $\zeta:\HFl(\A_H^F)\pf\HFl(\A_0^F)$ is well-defined via the homotopy $\beta^-_\infty(s)H$. Applying a homotopy of homotopies $\beta_r(s)H$ as in the proof of Theorem A shows that 
\beq
\zeta\circ\theta=\id:\HFl(\A_0^F)\pf\HFl(\A_0^F)\;,
\eeq
namely no breaking along non-trivial Reeb occurs during the homotopy. Hence $\theta$ is injective.
\end{proof}

\appendix

\section{$\A_H^F$ is generically Morse}\label{appendix}

In this appendix we prove Theorem \ref{thm:Morse}.

\subsection{Preparations}
The proof of the genericity of the Morse property follows a standard scheme, that is, once it is shown that a certain linear operator is surjective the theorem follows from Sard-Smale's theorem. Unfortunately, the standard approach by linearizing the functional using some connection leads to finding solutions of a rather complicated ODE on the manifold $M$. To circumvent this we first transform the problem and then in the end obtain a linear ODE in a vector space.

First, let us recall the definition of the perturbed Rabinowitz action functional
\bea
\A_H^F:\L\times\R&\pf\R\\
(v,\eta)&\mapsto-\int_0^1\lambda(v(t))[\partial_tv]-\int_0^1H(t,v)dt-\eta\int_0^1F(t,v)dt
\eea
where from now on $\L\equiv W^{1,2}(S^1,M)$ is the (completed) loop space of $M$. For convenience we abbreviate
\bea
\F:\L&\pf\R\\
v&\mapsto\int_0^1F(t,v)dt\;.
\eea
and
\beq
\A_H(v):=-\int_0^1\lambda(v(t))[\partial_tv]-\int_0^1H(t,v)dt\;.
\eeq
Thus, $\A_H^F(v,\eta)=\A_H(v)-\eta\F(v)$. We note that $\A^F_H(v,\eta)=\A_{\eta_0F+H}(v)+(\eta_0-\eta)\F(v)$, and therefore
\beq
d\A_H^F(v,\eta)[\hat{v},\hat{\eta}]=d\A_{\eta_0F+H}(v)[\hat{v}]-\hat{\eta}\F(v)+(\eta_0-\eta)d\F(v)[\hat{v}]
\eeq
where $\hat{v}\in\Gamma^{1,2}(v^*TM)$, the space of $W^{1,2}$ vector fields along $v$, and $\hat{\eta}\in\R$.
Hence at a critical point $x_0=(v_0,\eta_0)\in\Crit\A_H^F$ the Hessian equals
\beq
\He_{\A_H^F}(x_0)[(\hat{v}_1,\hat{\eta}_1),(\hat{v}_2,\hat{\eta}_2)]=\He_{\A_{\eta_0F+H}}(v_0)[\hat{v}_1,\hat{v}_2]-\hat{\eta}_1d\F(v_0)[\hat{v}_2]-\hat{\eta}_2d\F(v_0)[\hat{v}_1]\;.
\eeq
For a function $P:[0,1]\times M\pf\R$ and corresponding $\phi_P^1\in\Ham(M,\om)$ we define 
\beq\label{eqn:definition of L_H}
\L_P:=\{w\in W^{1,2}([0,1],M)\mid w(0)=\phi_P^1(w(1))\}\,,
\eeq
the twisted loop space, and introduce the diffeomorphism $\Phi_P:\L_P\pf\L$
\beq
\Phi_P(w)(t)=\phi_P^t(w(t))\;.
\eeq
For a fixed critical point $x_0=(v_0,\eta_0)$ of $\A_H^F$ we use this diffeomorphism to pull back $\A_H^F$ 
\beq
\At_{\eta_0,H}^F=(\Phi_{\eta_0 F+H}\times\id_\R)^*\A_H^F:\L_{\eta_0 F+H}\times\R\pf\R\;.
\eeq
We set $w_0:=\Phi_{\eta_0 F+H}^{-1}\circ v_0$, thus $w_0=$const. Then using $\big(\Phi_H^*d\A_H\big)(w)[\hat{w}]=\int \om(\partial_t w,\hat{w})$ we obtain
\beq
\He_{\At_{\eta_0,H}^F}(w_0,\eta_0)[(\hat{w}_1,\hat{\eta}_1),(\hat{w}_2,\hat{\eta}_2)]=\int_0^1\om(\partial_t\hat{w}_1,\hat{w}_2)dt-\hat{\eta}_1d\Ft(w_0)[\hat{w}_2]-\hat{\eta}_2d\Ft(w_0)[\hat{w}_1]\;.
\eeq
where $\Ft=\F\circ\Phi_{\eta_0 F+H}$. Using the special form of $F$ (see equation \eqref{eqn:def_of_F}) and $H\in \mathcal{H} $ (see equation \eqref{eqn:def_of_fancy_H}) we compute
\bea
\Ft(w)&=\int_0^1F(t,\phi^t_{\eta_0F+H}(w))dt=\int_0^{\tfrac12}F(t,\phi^t_{\eta_0F+H}(w))dt\\
&=\int_0^{\tfrac12}F(t,\phi^t_{\eta_0F}(w))dt=\int_0^{\tfrac12}F(t,w)dt\\
&=\int_0^1F(t,w)dt
\eea
Thus, the Hessian of $\At_{\eta_0,H}^F$ simplifies as follows (after integrating by parts)

\bea\label{eqn:Hessian_At_eta,H^F}
\He_{\At_{\eta_0,H}^F}&(w_0,\eta_0)[(\hat{w}_1,\hat{\eta}_1),(\hat{w}_2,\hat{\eta}_2)]\\
&=\int_0^1\om(\partial_t\hat{w}_1,\hat{w}_2)dt-\hat{\eta}_1\int_0^1dF(t,w_0)[\hat{w}_2]-\hat{\eta}_2\int_0^1dF(t,w_0)[\hat{w}_1]\\[1ex]
\eea

\subsection{The linearized operator}

We denote by $\mathcal{H}^k=\{H\in C^k(S^1\times M)\mid H(t,\cdot)=0\;\forall t\in[0,\tfrac12]\}$. Moreover, for $w\in\L_H$ (see equation \eqref{eqn:definition of L_H} for the definition) we define the bundle $\E_H\pf\L_H$ by
\beq
(\E_H)_w:=L^2([0,1],w^*TM)\;.
\eeq

\begin{Def}
Let $(v_0,\eta_0)$ be a critical point of $\A_H^F$ and $(w_0,\eta_0)$ the corresponding critical point of $\At_{\eta_0,H}^F$, that is, defined by the equation $v_0=\Phi_{\eta_0F+H}(w_0)$. Then we define the linear operator 
\beq
L_{(w_0,\eta_0,H)}:(T_{w_0}\L_{\eta_0F+H})\times\R\times\mathcal{H}\pf (\E_{\eta_0F+H})^\vee\times\R
\eeq
via the pairing with $(\hat{w}_2,\hat{\eta}_2)\in (\E_{\eta_0F+H})\times\R$
\bea
\langle L_{(w_0,\eta_0,H)}[\hat{w}_1,\hat{\eta}_1,\hat{H}], (\hat{w}_2,\hat{\eta}_2)\rangle&:=\He_{\At_{\eta_0,H}^F}(w_0,\eta_0)[(\hat{w}_1,\hat{\eta}_1),(\hat{w}_2,\hat{\eta}_2)]\\
&\quad+\int_0^1d((\Phi_{\eta_0F+H})^*\hat{H})(t,w_0)[\hat{w}_2(t)]dt
\eea
\end{Def}

\begin{Prop}\label{prop:linearized_operator_surjective}
The operator $L_{(w_0,\eta_0,H)}$ is surjective. In fact, $L_{(w_0,\eta_0,H)}$ is surjective when restricted to the space
\beq
\mathcal{V}:=\{(\hat{w},\hat{\eta},\hat{H})\in(T_{w_0}\L_{\eta_0F+H})\times\R\times\mathcal{H}\mid \hat{w}(\tfrac12)=0\}\;.
\eeq
\end{Prop}

\begin{Rmk}
The additional assertion of the surjectivity of $L_{(w_0,\eta_0,H)}|_{\mathcal{V}}$ is not used in the current article but will prove useful in the future. Since it added only two lines we decided to include it here.
\end{Rmk}

\begin{proof}
The $L^2$-Hessian is a self-adjoint Fredholm operator. Thus, the operator $L_{(w_0,\eta_0,H)}$ has closed image. Therefore, it suffices to prove that the annihilator of the image of $L_{(w_0,\eta_0,H)}$ vanishes. Let $(\hat{w}_2,\hat{\eta}_2)$ be in the annihilator of the image of $L_{(w_0,\eta_0,H)}$, that is 
\beq
\langle L_{(w_0,\eta_0,H)}[\hat{w}_1,\hat{\eta}_1,\hat{H}], (\hat{w}_2,\hat{\eta}_2)\rangle=0
\eeq
for all $(\hat{w}_1,\hat{\eta}_1,\hat{H})\in (T_{w_0}\L_{\eta_0F+H})\times\R\times\mathcal{H}$. This is equivalent to the following two equations:
\beq\label{eqn:A1}
\He_{\At_{\eta_0,H}^F}(w_0,\eta_0)[(\hat{w}_1,\hat{\eta}_1),(\hat{w}_2,\hat{\eta}_2)]=0\qquad\forall (\hat{w}_1,\hat{\eta}_1)\in (T_{w_0}\L_{\eta_0F+H})\times\R
\eeq
and
\beq\label{eqn:A2}
\int_0^1d\hat{H}_t(\phi_{\eta_0F+H}^t(w_0))[d\phi_{\eta_0F+H}^t(w_0)[\hat{w}_2]]=0 \qquad\forall \hat{H}\in\mathcal{H}
\eeq
Since the Hessian $\He_{\At_{\eta_0,H}^F}$ is a self-adjoint operator, equations  \eqref{eqn:Hessian_At_eta,H^F} and \eqref{eqn:A1} imply by elliptic regularity that $\hat{w}_2\in C^{k+1}([0,1],M)$ and satisfies the equation 
\beq\label{eqn:A3}
\partial_t\hat{w}_2-\hat{\eta}_2 X_F(t,w_0)=0
\eeq
and the linearized boundary condition
\beq
\label{eqn:A4}
\hat{w}_2(0)=d\phi_{\eta_0F+H}^1(w_0)[\hat{w}_2(1)]\;.
\eeq
In fact, when the Hessian is restricted to $\mathcal{V}$ then equation \eqref{eqn:A3} holds for all $t\neq\tfrac12$, since the Hessian is a local operator. Thus, by continuity, equation \eqref{eqn:A3} holds for all $t$ in any case.

From equation \eqref{eqn:A2} we deduce that
\beq\label{eqn:B}
\hat{w}_2(t)=0\quad\forall t\in[\tfrac12,1]\;.
\eeq
Using $F(t,x)=\rho(t)G(x)$ we rewrite equation \eqref{eqn:A3}
\beq
\partial_t\hat{w}_2-\hat{\eta}_2 \rho(t) X_G(w_0)=0\;.
\eeq
This is a linear ODE in the vector space $T_{w_0}M$ which we can solve
\beq\label{eqn:C}
\hat{w}_2(t)= \hat{w}_2(0)+\hat{\eta}_2\left(\int_0^t\rho(\tau)dt\right) X_G(w_0)\;.
\eeq
We recall (see equation \eqref{eqn:def_of_F}) that $F(t,x)=\rho(t)G(x)$ where $\int_0^t\rho(t)dt =1$ for all $t\in[\tfrac12,1]$. Combining this with equation \eqref{eqn:B} we conclude for $t\geq\tfrac12$
\beq\label{eqn:D}
0=\hat{w}_2(t)= \hat{w}_2(0)+\hat{\eta}_2 X_G(w_0)\;.
\eeq
Combining equations \eqref{eqn:A4} and \eqref{eqn:B} at $t=1$ we derive $\hat{w}_2(0)=0$. Hence, by equation \eqref{eqn:D} we have
\beq
\hat{\eta}_2 X_G(w_0)=0\;.
\eeq
Since $(w_0,\eta_0)$ comes from a critical point $(v_0,\eta_0)$ of $\A_H^F$ we know $G(v(0))=G(w_0)=0$, and therefore, $X_G(w_0)\neq0$ since $0$ was assumed to be a regular of $G$. In particular, 
\beq\label{eqn:E}
\hat{\eta}_2 =0
\eeq
Equations \eqref{eqn:C} and \eqref{eqn:E} immediately imply
\beq
\hat{w}_2(t)=0\quad\forall t\in[0,1]\;.
\eeq
Therefore, the annihilator of the image of $L_{w_0,\eta_0,H}$ vanishes and thus $L_{w_0,\eta_0,H}$ is surjective.
\end{proof}

\subsection{Proof of Theorem \ref{thm:Morse}}

We recall that $\L=W^{1,2}(S^1,M)$ and $\mathcal{H}^k=\{H\in C^k(S^1\times M)\mid H(t,\cdot)=0\;\forall t\in[0,\tfrac12]\}$. We define the Banach space bundle $\E\pf\L$ by $\E_v=L^2(S^1,v^*TM)$. We consider the section $S:\L\times\R\times\mathcal{H}^k\pf\E^\vee\times\R$ given by the differential of the Rabinowitz action functional $\A_H^F$ 
\beq
S(v,\eta,H):=d\A_H^F(v,\eta)\;.
\eeq 
where the perturbation $H\in\mathcal{H}^k$ is considered an additional variable. Its vertical differential $DS:T_{(v_0,\eta_0,H)}\L\times\R\times\mathcal{H}^k\pf\E_{(v_0,\eta_0,H)}^\vee$ at $(v_0,\eta_0,H)\in S^{-1}(0)$ is
\beq
DS_{(v_0,\eta_0,H)}[(\hat{v},\hat{\eta},\hat{H})]=\He_{\A_H^F}(v_0,\eta_0)\big[(\hat{v},\hat{\eta},\hat{H})\,;\:\bullet\:\big]+\int_0^1\hat{H}(t,v_0)dt
\eeq
Since the pull-back of $DS$ under the diffeomorphism $\Phi_{\eta_0F+H}\times\id_\R\times\id_{\mathcal{H}^k}$ is the operator $L_{(w_0,\eta_0,H)}$ in Proposition \ref{prop:linearized_operator_surjective}, the operator $DS$ is surjective. Thus, by the implicit function theorem the universal moduli space 
\beq
\M:=S^{-1}(0)
\eeq
is a smooth Banach manifold. We consider the projection $\Pi:\M\pf\mathcal{H}^k$. Then the $\A_H^F$ is Morse if and only if $H$ is a regular value of $\Pi$, which by the theorem of Sard-Smale form a generic set (for $k$ large enough). Moreover, the Morse condition is $C^k$-open. Thus, for functions in an open and dense subset of $\mathcal{H}^k$ the Rabinowitz action functional is Morse. Taking the intersection of all $k$ concludes the proof of Theorem \ref{thm:Morse}.\qed

%
%
\bibliographystyle{amsalpha}
\bibliography{../../../Bibtex/bibtex_paper_list}
\end{document}